\documentclass{amsart}

\usepackage{amsmath}
\usepackage{amscd}
\usepackage{amsfonts}
\usepackage{amssymb}
\usepackage{latexsym}
\usepackage{pifont}
\usepackage{eepic}
\usepackage{array}

\newcommand{\ncm}{\newcommand}

\ncm{\beq}{\begin{equation}}
\ncm{\eeq}{\end{equation}}
\ncm{\bea}{\begin{eqnarray}}
\ncm{\eea}{\end{eqnarray}}
\ncm{\beanon}{\begin{eqnarray*}}
\ncm{\eeanon}{\end{eqnarray*}}

\newtheorem{thm}{Theorem}[section]
\newtheorem{pro}[thm]{Proposition}
\newtheorem{lem}[thm]{Lemma}
\newtheorem{cor}[thm]{Corollary}

\theoremstyle{definition}
\newtheorem{defi}[thm]{Definition}

\theoremstyle{remark}

\numberwithin{equation}{section}

\def\MonCat{\mathsf{MonCat}}

\def\M{\mathsf{M}}
\ncm{\Cgd}{\mathsf{Cgd}}
\ncm{\BGD}{\mathsf{Bgd}}
\ncm{\BgdMap}{\mathsf{BgdMap}}
\ncm{\Bim}{\mathbb{BIM}}

\def\C{\mathcal{C}}

\ncm{\F}{\mathcal{F}}
\ncm{\G}{\mathcal{G}}

\ncm{\V}{\mathcal{V}}
\ncm{\W}{\mathcal{W}}

\ncm{\asso}{\mathbf{a}}
\ncm{\luni}{\mathbf{l}}
\ncm{\runi}{\mathbf{r}}

\ncm{\End}{\operatorname{End}}
\def\Hom{\mbox{\rm Hom}\,}

\ncm{\Mat}{\operatorname{Mat}}
\ncm{\Diag}{\operatorname{Diag}}
\ncm{\Ind}{\operatorname{Ind}}

\newcommand{\ci}{\circ}

\def\o{\otimes}
\def\x{\times}
\ncm{\amalgo}[1]{\underset{\scriptscriptstyle #1}{\o}}
\ncm{\mash}[1]{\,\raise2pt\hbox{$\underset{\scriptscriptstyle 
#1}{\scriptstyle\sharp}$}\,}

\ncm{\oA}{\amalgo{A}}
\ncm{\oB}{\amalgo{B}}
\ncm{\oR}{\amalgo{R}}
\ncm{\oT}{\amalgo{T}}
\ncm{\oS}{\amalgo{S}}
\ncm{\oN}{\amalgo{N}}
\ncm{\ex}[1]{\underset{\scriptscriptstyle #1}{\x}}
\def\bra{\langle}
\def\ket{\rangle}
\ncm{\rarr}[1]{\stackrel{#1}{\longrightarrow}}
\ncm{\larr}[1]{\stackrel{#1}{\longleftarrow}}
\def\cop{\Delta}

\def\eps{\varepsilon}
\ncm{\op}{\mathrm{op}}
\ncm{\coop}{\mathrm{coop}}
\ncm{\Fi}{\varphi}
\ncm{\IN}{\mathrm{in}}
\ncm{\OUT}{\mathrm{out}}
\ncm{\OR}{\overrightarrow}
\ncm{\OL}{\overleftarrow}

\def\oneT{^{(1)}}
\def\twoT{^{(2)}}
\def\threeT{^{(3)}}

\ncm{\I}{\mathcal{I}}

\def\du1{\hat 1}
\def\iso{\rarr{\sim}}
\def\ract{\triangleleft}
\def\lact{\triangleright}

\ncm{\under}{\mbox{\rm\_}\,}
\ncm{\Cnt}{\mathsf{C}}

\ncm{\ZZ}{\mathbb{Z}}
\ncm{\NN}{\mathbb{N}}

\ncm{\Oo}{\mathcal{O}}
\ncm{\Ha}{\mathcal{H}}
\ncm{\Ve}{\mathcal{V}}

\ncm{\sqrM}{$\sqrt{\text{Morita}}$\ }
\ncm{\ld}{\,.\,}
\ncm{\ud}{\,^.\,}
\ncm{\into}{\hookrightarrow}

\begin{document}

\title{Monoidal Morita equivalence}
\author[K. Szlach\'anyi]{Korn\'el Szlach\'anyi}
\address{Research Institute for Particle and Nuclear Physics, Budapest}
\email{szlach@rmki.kfki.hu}

\maketitle

Let $A$ be an algebra over the commutative ring $k$. 
It is well known that the category $\M_A$ of right $A$-modules is cocomplete, 
Abelian and the right regular object $A_A$ is a small projective generator. The 
latter three properties means precisely that  the functor 
$\Hom_A(A,\under):\M_A\to\M_k$ preserves coproducts, preserves cokernels and it 
is faithful, respectively. In fact this functor is monadic and has a right 
adjoint. It is also well known \cite{Mitchell} that the above properties 
characterize such categories: For a $k$-linear category $\C$ to be equivalent to 
$\M_A$ for some $k$-algebra $A$ it is sufficient that $\C$ is cocomplete, 
Abelian and possesses a small projective generator. Of course the algebra $A$ is 
determined by $\C$ only up to Morita equivalence. The analogue question for 
monoidal module categories has been studied by B. Pareigis in \cite{Pareigis: 
Morita}. With the advent of quantum groupoids it is worth reconsidering 
the question.

Therefore we are interested in monoidal structures on $\M_A$ admitting a strong 
monoidal forgetful functor to the category $_R\M_R$ of bimodules over some other 
$k$-algebra $R$. In the special case of $R=k$ one obtains that $A$ is a 
bialgebra \cite{Pareigis: nature}. The general case leads to bialgebroids 
\cite{Schauenburg: Bial}. The importance of module categories of bialgebroids, 
as well as those of Hopf algebroids, weak bialgebras and weak Hopf algebras, 
is that they provide examples of non-Tannakian monoidal categories. 

In a recent paper \cite{Sz: EM} the forgetful functors 
of bialgebroids over $R$ have been characterized abstractly as the strong 
monoidal monadic functors to $_R\M_R$ with a right adjoint. The bicategory of 
bialgebroids proposed in \cite{Sz: EM} is based on a comparison of certain 
monads on the base categories $_R\M_R$. Unfortunately, it does not admit a 
comparison of the monoidal categories $\M_A$ directly, without knowing at 
least a functor between their base categories. In Section \ref{sec: Bgd} a new 
construction of a bicategory of bialgebroids will be given which is completely 
base independent and in which the equivalence of two objects can be interpreted 
as monoidal Morita equivalence. The basic idea can be found already in 
Takeuchi's paper \cite{Takeuchi: sqrM}. This is why our monoidal Morita 
equivalence reduces to \sqrM equivalence in case of the Sweedler bialgebroids 
$E(R)$ which are the bialgebroid structures on the Sweedler corings $R^e=R^\op\o 
R$. Note that the monoidal version of classical Morita theory, what we study 
here, is \textit{a priori} different from (weak or not weak) Morita theory for 
monoidal categories \cite{Mueger: I}.

The results of Section \ref{sec: MM} put forward the notion of  strong comonoid 
progenerators: They classify the members of a monoidal  Morita equivalence class 
and their existence characterize the categories of modules over a bialgebroid. 
Similarly, existence of strong Frobenius progenerators characterize the 
categories $\M_A$ with $A$ being a Frobenius Hopf algebroid.

\section{The bicategory of coalgebroids}

We fix a commutative ring $k$ throughout the paper. All categories and functors 
are assumed to be $k$-linear even if not stated explicitly.

In his paper \cite{Takeuchi: sqrM} M. Takeuchi introduced the notion of an 
$R|S$-coring that play the role of 1-cells in a bicategory which we are going to
describe.
\begin{defi}
For $k$-algebras $R$ and $S$ an $R|S$-coalgebroid is an $R\o S$-$R\o S$-bimodule 
$C$ with an $S$-coring structure $\bra C,\cop,\eps\ket$ such that $\cop $ is 
also an $R$-$R$-bimodule map and its image is in the $R$-centralizer $(C\oS 
C)^R$. That is to say, using the familiar notation $\cop(c)=c\oneT\o c\twoT$,
\begin{align}
\cop(r\ld c\ld r')&=c\oneT\ld r'\oS r\ld c\twoT \label{cgd1}\\
r\ld c\oneT\oS c\twoT&=c\oneT\oS c\twoT\ld r \label{cgd2}
\end{align}
for all $c\in C$ and $r,r'\in R$.
\end{defi}
We denote the $S$-actions by upper dots and the 
$R$-actions by lower dots which becomes important in the special case when 
$R=S$.
In this case the ambiguous $(C\oS C)^S$ is denoted by $C\ex{S}C$.
 It follows from (\ref{cgd2}) that
\begin{equation}
\eps(r\ld c)=\eps(c\ld r)\qquad r\in R,\ c\in C\,.
\end{equation}
\begin{defi}
Let $\bra P,\cop_P,\eps_P\ket$ and $\bra Q,\cop_Q,\eps_Q\ket$ be 
$R|S$-coalgebroids. A $k$-linear map $\alpha\colon P\to Q$ satisfying
\begin{align*}
\alpha(s\ud p\ud s')&=s\ud\alpha(p)\ud s'\\
\alpha(r\ld p\ld r')&=r\ld\alpha(p)\ld r'\\
(\alpha\oS\alpha)\circ\cop_P&=\cop_Q\circ\alpha\\
\eps_Q\circ\alpha&=\eps_P
\end{align*}
is called a map of coalgebroids. The $R|S$-coalgebroids and $R|S$-coalgebroid 
maps form a category $\Cgd(S,R)$.
\end{defi}
\begin{defi}
Let $\Cgd$ be the bicategory with
\begin{enumerate}
\item objects being the $k$-algebras $R$, $S$, \dots,
\item hom-categories $\Cgd(S,R)$, 
\item horizontal composition $\Cgd(S,R)\x\Cgd(T,S)\to \Cgd(T,R)$ being given 
\begin{itemize}
\item for objects by $\bra P,Q\ket\mapsto P\odot Q$ where the $R|T$-coalgebroid 
$P\odot Q$ has underlying multimodule
\begin{align*}
P\odot Q&=P\o_{S^e} Q\\
&\equiv P\o Q/\{s\ud p\ud s'\o q-p\o s'\ld q\ld s\}\\
t\ud(p\odot q)\ud t'&=p\odot (t\ud q\ud t')\\
r\ld(p\odot q)\ld r'&=(r\ld p\ld r')\odot q
\end{align*}
and comultiplication and counit
\begin{align*}
\cop_{P\odot Q}(p\odot q)&=p\oneT\odot q\oneT\oT p\twoT\odot q\twoT\\
\eps_{P\odot Q}(p\odot q)&=\eps_Q(\eps_P(p)\ld q)
\end{align*}
\item and for arrows $\alpha:P\to P':S\to R$ and $\beta:Q\to Q':T\to S$ by
\begin{equation*}
(\alpha\odot\beta)(p\odot q)=\alpha(p)\odot\beta(q)
\end{equation*}
\end{itemize}
\item and with horizontal unit $E(R)$ at the object $R$ being the Sweedler 
coring $R^e$ as the $R|R$-coring with underlying multimodule
\begin{align*}
E(R)&=R\o R\\
r\ud(r_1\o r_2)\ud r'&=rr_1\o r_2r'\\
r\ld(r_1\o r_2)\ld r'&=r_1 r'\o rr_2
\end{align*}
and coring structure
\begin{align*}
\cop_{E(R)}(r_1\o r_2)&=(r_1\o 1_R)\oR(1_R\o r_2)\\
\eps_{E(R)}(r_1\o r_2)&=r_1r_2\,.
\end{align*}
\end{enumerate}
\end{defi}
The coherence isomorphisms
\[
O\odot(P\odot Q)\iso(O\odot P)\odot Q\quad \text{for}\ 
U\rarr{Q}T\rarr{P}S\rarr{O}R
\]
and
\[
E(R)\odot P\iso P,\quad P\odot E(S)\iso P\qquad\text{for}\ S\rarr{P}R
\]
are not written out explicitly. It suffices to note that they are the "same" as 
the coherence isomorphisms in the bicategory $\Bim(\M_k)$ of bimodules after 
forgetting the coring structure and considering the 1-cells $P\colon S\to R$ 
merely as $R^e$-$S^e$-bimodules. This can be formulated by saying that
the forgetful functor
\[
\mathbf{U}:\Cgd\to \Bim(\M_k)
\]
sending $\bra \,_{R^e}P_{S^e},\cop_P,\eps_P\ket$ to $_{R^e}P_{S^e}$ is a 
strict homomorphism of bicategories.

\section{The bicategory $\BGD$}  \label{sec: Bgd}

Bialgebroids have several, sligtly different though equivalent, formulations in 
the literature \cite{Takeuchi: x, Lu, Sz: Siena, Brz-Mil, Day-Street}.
For our purposes the following definition is the most appropriate: A bialgebroid
over $R$ is a monoid $\bra A,\mu,\eta\ket$ in the monoidal category $\Cgd(R,R)$. 
Especially, the monoidal unit $E(R)$ of 
$\Cgd(R,R)$ considered with the coherence isomorphism $E(R)\odot E(R)\iso E(R)$ 
as multiplication, i.e., with\[
(r_1\o r_2)\o_{R^e}(r'_1\o r'_2)\mapsto (r'_1r_1\o r_2r'_2)
\]
is a bialgebroid in this sense. 
So a natural way to introduce 1-cells $P:A\to B$ between bialgebroids and 
2-cells $\alpha:P\to Q:A\to B$ between such 1-cells is provided by considering 
the bicategory of bimodules in $\Cgd$.

\begin{defi}
The bicategory $\BGD$ of bialgebroids is defined as $\Bim(\Cgd)$. This means 
that the objects of $\BGD$ are the monoids\footnote{We mean monads, 
i.e., endo-1-cells with monoid structure, and not pseudomonoids which would be 
structures on 0-cells.}  in $\Cgd$,  the 1-cells $P:A\to B$ are the 
$B$-$A$-bimodules in $\Cgd$ and the 2-cells $\alpha:P\to Q:A\to B$ are the 
$B$-$A$-bimodule morphisms in $\Cgd$.
Horizontal composition of $B\rarr{P}A$ and $C\rarr{Q}B$ is denoted by $P\oB Q$.
\end{defi}
In the above definition we 
described the cells of $\BGD$ in the language of $\Cgd$. The more elementary 
language is that of $\Bim(\M_k)$ which we now use to describe $\BGD$.

\subsection{Bialgebroids}
A bialgebroid $A$ over $R$ is a monoid $\bra A,\mu,\eta\ket$ in $\Cgd(R,R)$. 
Therefore $A$ is a multimodule over $R$ with four actions denoted $r\ld a$, 
$r\ud a$, $a\ld r$ and $a\ud r$, respectively, and has a coring structure
\begin{align*}
\cop_A:&A\to A\oR A\\
\eps_A:&A\to R
\end{align*}
where $\oR$ is taken w.r.t. the upper dot actions. In order for $\bra \,
:\!\! A\colon,\cop,\eps\ket$ to be a coalgebroid the structure maps should 
satisfy\begin{align*}
\cop_A(r\ld a\ld r')&=a\oneT\ld r'\oR r\ld a\twoT \\
r\ld a\oneT\oR a\twoT&=a\oneT\oR a\twoT\ld r 
\end{align*}
for all $a\in A$ and $r,r'\in R$. 

The forgetting functor $\mathbf{U}$ applied to the monoid $\bra A,\mu,\eta\ket$ 
gives an $R^e$-ring. This means that $A$ is a $k$-algebra with 
multiplication $\bra a,a'\ket\mapsto aa'$ and unit element $1_A$ and the $\eta$
is a $k$-algebra map $R^e\to A$. So we can write
\[
\eta=t_A\o s_A:R^\op\o R\to A
\]
with uniquely determined algebra maps $s_A:R\to A$ and $t_A:R^\op\to A$ of 
commuting ranges. 

But $\eta$ is also a 2-cell $E(R)\to A$ in $\Cgd$ therefore it is a multimodule 
map
\begin{equation} \label{mmeta}
r\ud(1_A)\ud r'\ =t_A(r)s_A(r')\ =r'\ld(1_A)\ld r\qquad r,r'\in R
\end{equation}
and a coring map. Knowing already that $\cop_A$ and $\eps_A$ belong to $_R\M_R$, 
the latter reduces to the unitality conditions
\begin{equation}
\cop_A(1_A)=1_A\oR 1_A,\qquad\eps_A(1_A)=1_R\,.
\end{equation}

The $\mu$ being a 2-cell $A\odot A\to A$ it is a multimodule map
\begin{align}
(r\ld a\ld r')a'&=r\ld(aa')\ld r'\label{mmmul1}\\
a(r\ud a'\ud r')&=r\ud(aa')\ud r'\label{mmmul2}
\end{align}
and it is compatible with the coring structures therefore
\begin{align}
a\oneT{a'}\oneT\oR a\twoT{a'}\twoT&=(aa')\oneT\oR(aa')\twoT\\
\eps_A(aa')&=\eps_A(\eps_A(a)\ld a')\,.
\end{align}

Combining (\ref{mmmul1}) and (\ref{mmmul2}) with (\ref{mmeta})
 then imply 
\begin{align}
r\ud a\ud r'&=at_A(r)s_A(r')\label{mmfromst1}\\
r\ld a\ld r'&=s_A(r)t_A(r')a\label{mmfromst2}
\end{align}

Summarizing: A monoid $\bra A,\mu,\eta\ket$ in $\Cgd(R,R)$ is the same thing as 
a (right) bialgebroid $\bra A,R,s_A,t_A,\cop_A,\eps_A\ket$ satisfying the axioms 
given e.g. in \cite{Sz: Siena}. 

\subsection{Bialgebroid morphisms}

A 1-cell $P\colon B\to A$ from the bialgebroid $B$ over $S$ to the bialgebroid 
$A$ over $R$ has been defined as an $A$-$B$-bimodule in $\Cgd(S,R)$, i.e.,
a triple
\begin{align*}
P&\in \Cgd(S,R)\\
\lambda_P&:A\odot P\to P\\
\rho_P&:P\odot B\to P
\end{align*}
satisfying five commutative diagrams, as usual for bimodules.
Forgetting via $\mathbf{U}$ this bimodule becomes an ordinary bimodule $_AP_B$ 
in the (one object) bicategory $\M_k$. So $P$ is a $k$-module with left and 
right actions $\bra a,p\ket\mapsto a\lact p$ and $\bra p,b\ket\mapsto p\ract b$ 
of the $k$-algebras $A$ and $B$, respectively.
The multimodule structures on $A$ and $B$ allow to recognize the multimodule 
structure on $P$ because
\begin{align}
r\ld p\ld r'&=s_A(r)t_A(r')\lact p\\
s\ud p\ud s'&=p\ract t_B(s)s_B(s')\,.
\end{align}
The forgotten structures amount to have $S$-$S$ bimodule maps
\begin{align*}
\cop_P&:P\to P\oS P\\
\eps_P&:P\to S
\end{align*}
satisfying comonoid axioms in $_S\M_S$ and 
\begin{align}
r\ld p\oneT\oS p\twoT&=p\oneT\oS p\twoT\ld r \label{1-cell1}\\
\cop_P(a\lact p\ract b)&=a\oneT\lact p\oneT\ract b\oneT \oS a\twoT\lact 
p\twoT\ract b\twoT \label{1-cell2}\\
\eps_P(a\lact p\ract b)&=\eps_B(\eps_P(\eps_A(a)\ld p)\ld b) \label{1-cell3}
\end{align}
Note that  the last two equations express the fact that the $\lambda_P$ and 
$\rho_P$ are 2-cells in $\Cgd$.

Regrouping the axioms one obtains the following description of the 1-cells of 
$\BGD$.
\begin{lem} \label{lem: 1-cell}
The triple $\bra P,\cop_P,\eps_P\ket$ is a morphism of bialgebroids from $B$ 
over $S$ to $A$ over $R$ iff
\begin{enumerate}
\item $P$ is an $A$-$B$ bimodule $_AP_B$,
\item $\bra P_B,\cop_P,\eps_P\ket$ is a comonoid in $\bra\M_B,\oS,S_B\ket$,
\item $\cop_P:P\to(P\oS P)^R$ and
\item $\cop_P(a\lact p)=a\oneT\lact p\oneT\oS a\twoT\lact p\twoT$,
\item $\eps_P(a\lact p)=\eps_P(\eps_A(a)\ld p)$
\end{enumerate}
hold for all $a\in A$ and $p\in P$.
\end{lem}
Loosely speaking, 1-morphisms $B\to A$ of bialgebroids are right $B$-module 
coalgebras with compatible left $A$-module structure.

\subsection{Bialgebroid transformations}

If both $P$ and $Q$ are bialgebroid morphisms from $B$ over $S$ to $A$ over $R$ 
then a 2-cell $\alpha\colon P\to Q\colon B\to A$ is nothing but an 
$A$-$B$-bimodule map
\[
\alpha\colon\ _AP_B\to\ _AQ_B
\]
satisfying
\begin{align*}
\alpha(p\oneT)\oS\alpha(p\twoT)&=\alpha(p)\oneT\oS\alpha(p)\twoT\\
\eps_Q(\alpha(p))&=\eps_P(p)\,.
\end{align*}

\section{Embedding into $\MonCat$}

If $A$ is a right bialgebroid over $R$ then the category $\M_A$ of right 
$A$-modules has a unique monoidal structure $\bra \M_A,\oR,R_A\ket$ such that 
the forgetful functor $U_A:\M_A\to\,_R\M_R$ associated to the algebra map
$\eta:R^e\to A$ is strict monoidal. This is the object map of a morphism of 
bicategories.
\begin{defi}
Let $\Ha:\BGD\to \MonCat^\op$ be the morphism of bicategories into the 
2-category of monoidal categories, monoidal functors and monoidal natural 
transformations
\begin{enumerate}
\item which maps the bialgebroid $A$ to the monoidal category $\M_A$,
\item the 1-cell $P:B\to A$ to the monoidal functor 
\begin{align*}
\Ha(P)&:=\Hom_B(P,\under):\M_B\to\M_A\\
\Ha(P)_{M,N}&:\Hom_B(P,M)\oR\Hom_B(P,N)\to\Hom_B(P,M\oS N)\\
&\mu\oR\nu\mapsto(\mu\oS\nu)\circ\cop_P\\
\Ha(P)_0&:R_A\to\Hom_B(P,S)\\
&r\mapsto \eps_P(r\ld\under)
\end{align*}
\item and the 2-cell $\alpha:P\to Q:B\to A$ to the monoidal natural 
transformation
\[
\Hom_B(\alpha,M):\Hom_B(Q,M)\to\Hom_B(P,M)\,.
\]
\end{enumerate}
\end{defi}
For each $C\rarr{Q}B\rarr{P}A$ the natural isomorphism
\[
\Hom_B(P,\Hom_C(Q,M))\iso\Hom_C(P\oB Q,M)
\]
and for each $A$ the natural isomorphism
\[
\M_A\iso\Hom_A(A,\under)
\]
endow $\Ha$ with the structure of a homomorphism of bicategories 
\cite{Leinster: BB}.

\begin{lem} \label{H is ff}
$\Ha$ is locally faithful and full.
\end{lem}
\begin{proof}
Let $\alpha:P\to Q:B\to A$ be a 2-cell and assume that $\Ha(\alpha)=0$. Then
$\mu\circ\alpha=0$ for all $\mu\in\Hom_B(Q,M)$ and for all $M\in\M_B$. Choosing 
$M=Q$ and $\mu=Q$ we obtain $\alpha=0$.
Now let $\kappa:\Ha(Q)\to\Ha(P)$ be any monoidal natural transformation.
Then Yoneda Lemma  implies that there is an $A$-$B$-bimodule map 
$\alpha:P\to Q$ such that $\kappa_M=\Hom_B(\alpha,M)$. Since $\kappa$ is 
monoidal, 
\begin{align*}
(\alpha\oS\alpha)\circ\cop_P&=\Ha(P)_{Q,Q}\circ(\kappa_Q\oR\kappa_Q)(Q\oR Q)\\
&=\kappa_{Q\oS Q}\circ\Ha(Q)_{Q,Q}(Q\oR Q)=\cop_Q\circ\alpha
\end{align*}
and
\begin{equation*}
\eps_Q\circ\alpha=\kappa_S\circ\Ha(Q)_0(1_R)=\Ha(P)_0(1_R)=\eps_P\,.
\end{equation*}
Thus $\alpha$ is a map of coalgebroids. So we have proven that the functor
\[
\Ha(B,A)\colon\BGD(B,A)\to\MonCat^\op(\M_B,\M_A)
\]
is full and faithful for any pair of bialgebroids $A$, $B$.
\end{proof}
\begin{lem} \label{theH(P)s}
Let $A$ and $B$ be bialgebroids over $R$ and $S$, respectively, and let 
$G:\M_B\to\M_A$ be a monoidal functor. Then $G$ is isomorphic to $\Ha(P)$ for 
some bialgebroid morphism $P:B\to A$ if and only if the underlying ordinary 
functor of $G$ has a left adjoint.
\end{lem}
\begin{proof}
Clearly the functor $\Hom_B(P,\under)$ has left adjoint, namely $\under\oA 
P$.Now assume that $\bra G,G_2,G_0\ket$ is a monoidal functor and $F\dashv G$ is 
a left adjoint. Then $F$ has an opmonoidal structure $\bra F,F^2,F^0\ket$ and 
there is a monoidal natural isomorphism
\begin{equation*}
G(M)\cong \Hom_A(A,G(M))\cong\Hom_B(F(A),M)\equiv\Ha(P)(M)
\end{equation*}
where $P:=F(A)$ is an $A$-$B$-bimodule and has coring structure
\begin{align*}
\cop_P&=F^{A,A}\circ F(\cop_A)\\
\eps_P&=F^0\circ F(\eps_A)
\end{align*}
which satisfies the compatibility conditions (\ref{1-cell1}), (\ref{1-cell2}) 
and(\ref{1-cell3}).
\end{proof}
The characterization of the objects $\Ha(A)$ within $\MonCat$ is postponed until 
Section \ref{sec: MM}.

\section{Special morphisms}

\subsection{Bialgebroid maps}
Let $A$ and $B$ be bialgebroids over $R$ and $S$, respectively. A bialgebroid 
map $\bra f,f_0\ket\colon A\to B$ is a pair of algebra maps $f:A\to B$ and 
$f_0:R\to S$ such that 
\begin{align*}
f\circ s_A&=s_B\circ f_0\\
f\circ t_A&=t_B\circ f_0\\
(f\oS f)\circ\cop_A&=\cop_B\circ f\\
\eps_B\circ f&=f_0\circ\eps_A
\end{align*}
where notice that  $f_0$ is uniquely determined by $f$ via $f_0=\eps_B\circ 
f\circ s_A$. Therefore we shall often say that "$f$ is a bialgebroid map" 
without mentioning $f_0$.

The bialgebroids and bialgebroid maps form a category $\BgdMap$ and 
the notion of isomorphism in this category leads to the
\begin{defi}
Two bialgebroids $A$ over $R$ and $B$ over $S$ are called isomorphic if there
exist bialgebroid maps $f:A\to B$ and $g:B\to A$ such that $f\circ g$ and 
$g\circ f$ are identities.
\end{defi}

To any bialgebroid map $\bra f,f_0\ket\colon A\to B$ we can associate a 
bialgebroid morphism $f^*$ as follows. As an $A$-$B$-bimodule it is $B$ with 
left action of $A$ induced by $f$. The coring structure is 
inherited from $B$. It is easy to check that the triple
$\bra f^*=\,_{f(A)}B_B,\cop_{f^*}=\cop_B,\eps_{f^*}=\eps_B\ket$ satisfy the 
axioms for a 1-cell in $\BGD$. As a matter of fact, 
\begin{align*}
r\ld b\oneT\oS b\twoT&=f(s_A(r))b\oneT\oS b\twoT=s_B(f_0(r))b\oneT\oS b\twoT\\
&=b\oneT\oS t_B(f_0(r))b\twoT=b\oneT\oS f(t_A(r))b\twoT\\
&=b\oneT\oS b\twoT\ld r\\
\cop_{f^*}(a\lact b)&=\cop_B(f(a))\cop_B(b)=(f\oS f)(\cop_A(a))\cop_B(b)\\
&=a\oneT\lact b\oneT\oS a\twoT\lact 
b\twoT\\\eps_{f^*}(a\lact b)&=\eps_B(f(a)b)=\eps_B(\eps_B(f(a))\ld b)= 
\eps_B(f_0(\eps_A(a))\ld b)\\&=\eps_{f^*}(\eps_A(a)\ld b)
\end{align*}

\subsection{The forgetful functors}

For $A$ a bialgebroid over $R$ the algebra map $\eta:R^e\to A$ defines a 
forgetful functor $U_A:\M_A\to\M_{R^e}$ which, when $\M_{R^e}$ is identified 
with the monoidal category $_R\M_R$, is a strict monoidal functor.
One expects that this functor is the same as $\Ha(\eta^*)$ for the 1-cell 
$\eta^*: A\to E(R)$ associated to the bialgebroid map $\eta:E(R)\to A$.
\begin{pro}
Let $\bra A,\mu,\eta\ket$ be a bialgebroid over $R$. Then
\begin{enumerate}
\item $\eta:E(R)\to A$ is a bialgebroid map
\item and there is a monoidal isomorphism $U_A\cong\Ha(\eta^*)$.
\end{enumerate}
\end{pro}
\begin{proof}
(1) $\eta$ is an algebra map $R^e\to A$, $(r_1\o r_2)\mapsto t_A(r_1)s_A(r_2)$. 
Writing simply $E$ for $E(R)$ we have $s_E(r)=1_R\o r$ and $t_E=r\o 1_R$
therefore $\eta_0:=\eps_A\circ\eta\circ s_E\colon R\to R$ is the identity and we 
have
\begin{align*}
\eta\circ s_E&=s_A\\
\eta\circ t_E&=t_A\\
(\eta\oR\eta)\circ \cop_E&=\cop_A\\
\eps_A\circ\eta&=\eps_E
\end{align*}
(2) The 1-morphism $\eta^*\colon A\to E(R)$ has underlying bimodule $_{R^e}A_A$
with left action $(r_1\o r_2)\lact a=r_2\ld a\ld r_1$.
The functor $\Ha(\eta^*)$ maps the right $A$-module $M$ to $\Hom_A(A,M)$ 
therefore the right $R^e$-module maps
\begin{align*}
\nu_M&\colon U_A(M)\to\Ha(\eta^*)(M)\\
&m\mapsto\{a\mapsto m\ract a\}
\end{align*}
which are natural in $M$ define the required natural isomorphism. It remains to 
show that $\nu$ is monoidal, i.e.,
\begin{equation*}
\begin{CD}
U_A(M)\oR U_A(N)@=U_A(M\oR N)\\
@V{\nu_M\oR\nu_N}VV @VV{\nu_{M\oR N}}V\\
\Ha(\eta^*)(M)\oR\Ha(\eta^*)(N)@>>(\under\oR\under)\circ\cop_A>\Ha(\eta^*)(M\oR 
N)
\end{CD}
\end{equation*}
and
\begin{equation*}
\begin{CD}
_RR_R@=U_A(R)\\
@| @VV{\nu_R}V\\
_RR_R@>>\Ha(\eta^*)_0>\Ha(\eta^*)(R)
\end{CD}
\end{equation*}
are commutative in $_R\M_R$. 
The first diagram evaluated on $m\oR n$ is equivalent to $(m\oR n)\ract a=
(m\ract a\oneT)\oR(n\ract a\twoT)$ and the second one is equivalent to the 
formula $r\ract a=\eps_A(r\ld a)$ for the trivial representation.
\end{proof}

\subsection{When $\Ha(P)$ is strong monoidal}\label{ss: H(P) strong}

Let $P\colon B\to A$ be a morphism of bialgebroids and assume that 
$\Ha(P)\colon\M_B\to\M_A$ is strong monoidal, i.e., the $\Ha(P)_{M,N}$ for all 
objects $M$, $N$ and the $\Ha(P)_0$ are isomorphisms. Such situations
motivate the 
\begin{defi}
A comonoid $\bra g,\gamma,\pi\ket$ in a closed monoidal category $\bra 
\C,\boxtimes,e,\asso,\luni,\runi\ket$ is called a strong comonoid if the 
monoidal functor 
\begin{align*}
\Ha=\Hom(g,\under):\C&\to \,_T\C_T\\
\Ha_{a,b}:\Hom(g,a)\oT\Hom(g,b)&\to\Hom(g,a\boxtimes b),\quad 
\alpha\oT\beta\mapsto (\alpha\boxtimes\beta)\circ\gamma\\
\Ha_0:T&\to\Hom(g,e),\quad \tau\mapsto \tau
\end{align*}
is strong. Here $T$ is the convolution monoid $\Hom(g,e)$ and for each object 
$a$ the $\Hom(g,a)$ is given the $T$-$T$-bimodule structure
$\tau\cdot \alpha=\luni_a\circ(\tau\boxtimes \alpha)\circ\gamma$, $\alpha\cdot 
\tau=\runi_a\circ(\alpha\boxtimes \tau)\circ\gamma$.
\end{defi}
For $P$ a strong comonoid in $\M_B$ there is a bialgebroid 
$\bra E,T,s_E,t_E,\cop_E,\eps_E\ket$, called the endomorphism bialgebroid, 
defined as follows.
\begin{align}
E&:=\End_B(P)\quad\text{as an algebra}\label{E1}\\
T&:=\Hom_B(P,S)\quad \text{with multiplication } \tau\star \tau':=(\tau\oS 
\tau')\circ\cop_P\\
s_E&: \tau\mapsto(P\oS \tau)\circ\cop_P\\
t_E&\colon \tau\mapsto (\tau\oS P)\circ\cop_P\\
\cop_E&:\alpha\mapsto \Ha_{P,P}^{-1}(\cop_P\circ\alpha)\\
\eps_E&:\alpha\mapsto \eps_P\circ \alpha\,,\label{Elast}
\end{align}
If $P_B$ is the restriction of a $P:B\to A$ then the $T$ and $R$ actions on $E$ 
are related by the isomorphism $\Ha(P)_0\colon R\iso T$. Namely, for $r\in R$, 
$\tau=\eps_P(r\ld\under)$ we have
\begin{align*}
\tau\ud\alpha&=\alpha\circ t_E(\tau)=\alpha(\under\ld r)\equiv r\ud\alpha\\
\alpha\ud\tau&=\alpha\circ s_E(\tau)=\alpha(r\ld\under)\equiv \alpha\ud r\ .
\end{align*}
\begin{lem} \label{lem: strong H(P)}
If the bialgebroid morphism $P:B\to A$ is such that $\Ha(P)=\Hom_B(P,\under)$ is 
a strong monoidal functor then 
\begin{enumerate}
\item $E=\End_B(P)$ is a bialgebroid with the structure maps given above, 
\item with the natural action $\alpha\lact p=\alpha(p)$ of $E$ on $P$ the triple 
$\bra \,_EP_B,\cop_P,\eps_P\ket$ is a bialgebroid morphism $B\to E$ and
\item  the map $\lambda:A\to E$ given by $_AP$ is a map of 
bialgebroids.
\end{enumerate}
\end{lem}
\begin{proof}
(1) This follows by patiently substituting (\ref{E1})-(\ref{Elast}) into the 
(right) bialgebroid axioms of \cite{Sz: Siena}.

(2) $P$ is an $E$-$B$-bimodule by construction of $E$. The coring structure of 
$P$ is compatible with the $B$-action since $P$ is a bialgebroid morphism from 
$B$. So we are left with proving compatibility with the $E$-action, i.e., 
equations  (3), (4), (5) of Lemma \ref{lem: 1-cell}.
\begin{align*}
\tau\ld p\oneT\oS p\twoT&=p\oneT\ud \tau(p\twoT)\oS p\threeT=
p\oneT\oS \tau(p\twoT)\ud p\threeT\\
&=p\oneT\oS p\twoT\ld \tau\\
\alpha\oneT(p\oneT)\oS \alpha\twoT(p\twoT)&=(\alpha\oneT\oS\alpha\twoT)\circ 
\cop_P(p)\\
&=\cop_P(\alpha( p))\\
\eps_P(\alpha(p))&=\eps_E(\alpha)(p)=\eps_E(\alpha)(p\oneT)\eps_P(p\twoT)\\
&=\eps_P(\eps_E(\alpha)(p\oneT)\ud p\twoT)=\eps_P(\eps_E(\alpha)\ld p)\,.
\end{align*}

(3) The map $\lambda\colon a\mapsto a\lact\under$ is an algebra map 
obviously. Let $\lambda_0:=\eps_E\circ \lambda\circ s_A: R\to T$ then
$\lambda_0(r)=\eps_P(r\ld\under)=\Ha(P)_0(r)$, thus $\lambda_0=\Ha(P)_0$ is an 
isomorphism. Now it is clear that
\begin{align*}
\lambda\circ s_A(r)&=r\ld\under=s_E\circ\lambda_0(r)\\
\lambda\circ t_A(r)&=\under\ld r=t_E\circ\lambda_0(r)\,.
\end{align*}
Using compatibility of $\cop_P$ with the $A$-action, then  the restricted 
naturality of $\Ha(P)_{P,P}$  we obtain
\begin{align*}
\cop_E(\lambda(a))&=m^{-1}(\cop_P\circ\lambda(a))=m^{-1}\left((\lambda(a\oneT) 
\oS\lambda(a\twoT))\circ\cop_P\right)\\
&=\lambda(a\oneT)\oT\lambda(a\twoT)
\end{align*}
and finally
\begin{equation*}
\eps_E(\lambda(a))=\eps_P\circ\lambda(a)=\eps_P(\eps_A(a)\ld\under)=
\lambda_0(\eps_A(a))\,.
\end{equation*}
\end{proof}

\section{Monoidal Morita equivalence} \label{sec: MM}

In this section we study the properties of bialgebroid morphisms $P:B\to A$ that 
give rise to monoidal category equivalences $\Ha(P):\M_B\simeq\M_A$.
By Lemma \ref{theH(P)s} all monoidal equivalences are obtained in this way.

Forgetting the monoidal structure, at first, we are in the situation of 
classical Morita theory. The functor 
$\Hom_B(P,\under):\M_B\to\M_A$ is an equivalence of categories, so is its left 
adjoint $\under\oA P:\M_A\to\M_B$, therefore $_AP_B$ is a Morita equivalence 
bimodule with inverse equivalence bimodule $_BQ_A=\Hom_B(P,B)$. 
It follows that $_AP_B$ and $_BQ_A$ are 
faithfully balanced bimodules and the $P_B$, $_AP$, $Q_A$, $_BQ$ are all 
progenerators (i.e., finitely generated projective generators). In particular 
$A$ is determined up to isomorphism by $P_B$ as $\End_B(P)$. 

The monoidal structure on the functor $\Ha(P)$ imposes on the progenerator $P_B$ 
a comonoid structure $\bra P,\cop_P,\eps_P\ket$. Since a monoidal equivalence 
functor is nothing else but an equivalence functor with a strong monoidal 
structure, this comonoid must be strong and every strong comonoid progenerator
$\bra P,\cop_P,\eps_P\ket$ determines a monoidal equivalence 
$\Hom_B(P,\under):\M_B\to\M_E$ where $E$ is the endomorphism bialgebroid of 
$P$. 

\begin{defi}
Over the base category $\M_k$ of $k$-modules consider two bialgebroids: $A$ 
over $R$ and $B$ over $S$. We say that $A$ and $B$ are monoidally Morita 
equivalent and write $A\simeq B$ if there is a $k$-linear monoidal 
category equivalence $\M_A\simeq\M_B$.
\end{defi}
Monoidal Morita equivalence is the same as 
equivalence of objects in the bicategory $\BGD$. 
\begin{lem} 
The Morita equivalence class of a bialgebroid $B$ can be represented in the 
monoidal module category $\M_B$ as follows.
\begin{enumerate} 
\item Let $A\simeq B$ and let $P:B\to A$ be an equivalence in $\BGD$. Then $\bra 
P,\cop_P,\eps_P\ket$ is a strong comonoid progenerator in $\M_B$ 
and $\lambda_P:A\to\End(P_B)$ is an isomorphism of bialgebroids.
\item If $\bra P,\cop_P,\eps_P\ket$ is a strong comonoid progenerator in 
$\M_B$then $\End(P_B)\simeq B$.
\end{enumerate}
\end{lem}
\begin{proof}
(1) $P$ is an equivalence in $\BGD$ iff $\Ha(P)$ is a monoidal equivalence in 
$\MonCat$. In particular $\Ha(P)$ is strong monoidal. Thus Lemma \ref{lem: 
strong H(P)} implies that $P$ is a strong comonoid progenerator and $\lambda_P$ 
is an isomorphism.
(2) Since $P_B$ is a progenerator, $_EP_B$ is a Morita 
equivalence bimodule. Therefore $\Ha(P)$ is an equivalence of categories and is 
equipped with a strong monoidal structure. Therefore $\Ha(P):\M_B\to\M_E$ is a 
monoidal equivalence.
\end{proof}

The above characterization of the Morita 
equivalence class of $B$ would not be very useful without the next proposition 
which allows one to replace the global property of a comonoid being "strong"  
with a local one.
\begin{pro} \label{pro: strongness}
Let $F:\bra\C,\boxtimes,e\ket\to\bra\M,\o,i\ket$ be a colimit preserving 
additive monoidal functor between closed monoidal Abelian categories.
We assume also that $\C$ is cocomplete and possesses a generator $g$.
Then $F$ is strong monoidal iff the two arrows $F_{g,g}:Fg\o Fg\to F(g\boxtimes 
g)$ and $F_0:i\to Fe$ are isomorphisms in $\M$.
\end{pro}
\begin{proof}
Only the "if" part requires proof.
Every object $c$ of $\C$ has a presentation
\[
\coprod^J g\rarr{\alpha}\coprod^I g\rarr{\beta} c\to 0
\]
for some small sets $I$ and $J$. So the proof consists of two steps: At first we 
show that $F_{a,b}$ are invertible for $a$ and $b$ being arbitrary coproducts of 
$g$. At second we show invertibility of $F_{a,b}$ for $a$ and $b$ being 
cokernels of arrows like $\alpha$ above. 

(1) Let $u_i:g\to a$, $i\in I$  be injections of the coproduct $a$ and 
assume that invertibility of $F_{g,b}$ has already been proven. The arrows 
\[
\begin{CD}
Fg\o Fb@>Fu_i\o Fb>>Fa\o Fb
\end{CD}
\qquad i\in I
\]
are injections of a coproduct in $\M$. The universal property of this coproduct 
implies that $F_{a,b}$ is the unique arrow making the naturality diagram
\begin{equation}
\begin{CD}
Fg\o Fb@>F_{g,b}>>F(g\boxtimes b)\\
@V{Fu_i\o Fb}VV @VV{F(u_i\boxtimes b)}V\\
Fa\o Fb@>F_{a,b}>>F(a\boxtimes b)
\end{CD}
\end{equation}
commute. But also the arrows
\[
\begin{CD}
F(g\boxtimes b)@>F(u_i\boxtimes b)>>F(a\boxtimes b)
\end{CD}
\qquad i\in I
\]
are injections of a coproduct therefore there exists a unique $F'_{a,b}$ making
\begin{equation}
\begin{CD}
Fg\o Fb@<F_{g,b}^{-1}<<F(g\boxtimes b)\\
@V{Fu_i\o Fb}VV @VV{F(u_i\boxtimes b)}V\\
Fa\o Fb@<F'_{a,b}<<F(a\boxtimes b)
\end{CD}
\end{equation}
commutative. Now applying the universality properties of the two coproducts 
again we find that $F'_{a,b}\circ F_{a,b}$ and $F_{a,b}\circ F'_{a,b}$ are 
identities, thus $F_{a,b}$ is invertible. Using the above argument at first with 
$b=g$ we obtain that $F_{a,g}$ is invertible if $a\cong\coprod^I g$. Thus, 
interchanging the roles of the two tensorands, also $F_{g,b}$ is invertible if 
$b\cong\coprod^J g$. Then using the argument at the second time invertibility of 
$F_{a,b}$ follows for all $a$ and $b$ that are coproducts of $g$.

(2) Let $a\rarr{\alpha}b\rarr{\beta}c\to 0$ be exact where $a\cong\coprod^Ig$,
$b\cong\coprod^J g$ and assume that $d$ is such that invertibility of $F_{a,d}$,
$F_{b,d}$ is already proven. Then the columns of the diagram
\begin{equation}
\begin{CD}
Fa\o Fb@>F_{a,d}>>F(a\boxtimes d)\\
@V{F\alpha\o Fd}VV @VV{F(\alpha\boxtimes d)}V\\
Fb\o Fb@>F_{b,d}>>F(b\boxtimes d)\\
@V{F\beta\o Fd}VV @VV{F(\beta\boxtimes d)}V\\
Fc\o Fb@>F_{c,d}>>F(c\boxtimes d)\\
@VVV @VVV\\
0@= 0
\end{CD}
\end{equation}
are exact, too. Therefore $F_{c,d}$ is the unique arrow factorizing 
$F(\beta\boxtimes d)\circ F_{b,d}$ through the cokernel $F\beta\o Fd$. 
Similarly, there is a unique arrow factorizing $(F\beta\o Fd)\circ F_{b,d}^{-1}$ 
through the cokernel $F(\beta\boxtimes d)$, which is then necessarily the 
inverse of $F_{c,d}$.
Applying this argument at first for $d\cong \coprod^K g$ we obtain invertibility 
of $F_{c,d}$ and, interchanging the roles of the tensorands, of $F_{d,c}$, too.
Finally, we can apply the argument for arbitrary object $d$ which finishes the 
proof.
\end{proof}

\begin{cor} \label{cor: strongness}
Let $B$ be a bialgebroid over $S$ and let $\bra P,\cop_P,\eps_P\ket$ be a 
comonoid in $\M_B$.
\begin{enumerate}
\item If $P_B$ is finitely generated projective then $\bra 
P,\cop_P,\eps_P\ket$ is a strong comonoid iff 
\begin{align*}
Q\oT Q&\to\Hom_B(P,B\oS B)\\
q\oT q'&\mapsto (q\oS q')\circ\cop_P
\end{align*}
is an isomorphism, where $Q=\Hom_B(P,B)$, $T=\Hom_B(P,S)$ and the bimodule 
structure of $Q$ is given by $(t\ud q)(p)=t(p\oneT)\ud q(p\twoT)$, 
$(q\ud t)(p)=q(p\oneT)\ud t(p\twoT)$.
\item If $P_B$ is a progenerator then $\bra 
P,\cop_P,\eps_P\ket$ is a strong comonoid iff 
\begin{align*}
E\oT E:&\to \Hom_B(P,P\o_S P)\\
\alpha\oT\beta&\mapsto (\alpha\oS\beta)\circ\cop_P
\end{align*}
is an isomorphism, where $E=\End_B(P)$ is a
$T$-$T$-bimodule over the convolution monoid $T=\Hom_B(P,S)$ as in Subsection 
\ref{ss: H(P) strong}.
\end{enumerate}
\end{cor}

Now we have all the necessary tools to give a characterization of the categories 
of modules of a bialgebroid.
\begin{thm} \label{thm: main}
For a monoidal category $\C$ the following conditions are equivalent:
\begin{enumerate}
\item There is a bialgebroid $A$ and a monoidal equivalence $\C\simeq\M_A$.
\item There is an algebra $R$ and a strong monoidal monadic left adjoint functor 
$U:\C\to\,_R\M_R$.
\item $\C$ is cocomplete, Abelian with a small projective generator  
admitting a strong comonoid structure.
\item $\C$ is closed, cocomplete, Abelian with a small projective generator  
admitting a comonoid structure $\bra g,\gamma,\pi\ket$ such that
\begin{align*}
\End g\ \oT\ \End g&\to \Hom(g, g\boxtimes g)\\
\alpha\oT\beta&\mapsto (\alpha\boxtimes \beta)\ci\gamma
\end{align*}
is an isomorphism in $_T\M_T$ where $T=\Hom(g,e)$ with the convolution monoid
structure.
\end{enumerate}
\end{thm}

\begin{proof}
$(1)\Leftrightarrow (2)$ For a bialgebroid $A$ over $R$ the forgetful functor 
$U:\M_A\to\,_R\M_R$ is such a functor. If $R$ and $U$ is given then 
the bialgebroid $A$ can be reconstructed as the representing object of the 
opmonoidal monad $\under\o_{R^e}A$ on $_R\M_R$. For details see \cite{Sz: EM} 
for left bialgebroids though.

$(1)\Rightarrow (3)$  The $\M_A$ is cocomplete Abelian with $A_A$ a small 
projective generator. The coring structure $\bra A_A,\cop_A,\eps_A\ket$ gives 
the functor $\Hom_A(A,\under):\M_A\to \,_R\M_R$ a strong monoidal structure 
where $R$ is the convolution monoid $\Hom_A(A,R)$.

$(3)\Rightarrow (1)$ By Mitchell's Theorem the functor 
$\Hom(g,\under):\C\to\M_A$ where $A=\End g$ is an equivalence of categories 
if $g$ is a small projective generator. If $g$ has also a strong comonoid 
structure then the endomorphism bialgebroid construction of Subsection \ref{ss: 
H(P) strong} equips $A$ with right bialgebroid structure over $R=\Hom(g,e)$  
such that $\Hom(g,\under)$ is strong monoidal. 
Hence it is a monoidal equivalence $\C \simeq\M_A$.

$(4)\Rightarrow(3)$ Apply Proposition \ref{pro: strongness} to the normal 
monoidal functor $F=\Hom(g,\under):\C\to\,_T\M_T$ to conclude that $\bra 
g,\gamma,\pi\ket$ is strong.

$(3)\Rightarrow(4)$ Since already $(3)\Rightarrow(1)$ and $\M_A$ is closed 
\cite{Schauenburg: Duals},  so is $\C$.
\end{proof}

Let us briefly discuss the case of Frobenius Hopf algebroids.
An invariant for an object $g$ in a monoidal category $\bra\C,\boxtimes,e\ket$ 
is an arrow $\iota:e\to g$. A comonoid $\bra g,\gamma,\pi\ket$ is called 
Frobenius if $g$ is selfdual and admits a Frobenius integral. The latter means 
an invariant $\iota$ such that $\gamma\ci\iota$ 
is coevaluation for some evaluation $g\boxtimes g\to e$
exhibiting $g$ as its own left (and therefore right) dual.
This terminology is justified by the observation: Frobenius integrals exists for 
a comonoid if and only if there is an extension of the comonoid structure to a 
Frobenius algebra $\bra g,\nu,\iota,\gamma,\pi\ket$ in $\C$. 
If moreover the comonoid is strong then the endomorphism bialgebroid 
$E=\End g$ becomes equipped with antipode 
\begin{align*}
S_E(\alpha)&=\runi_g\ci(g\boxtimes\eps)\ci(g\boxtimes\nu)\ci(g\boxtimes
(\alpha\boxtimes g))\ci\asso_{g,g,g}^{-1}\ci\\
&(\gamma\boxtimes g)\ci(\iota\boxtimes g)\ci\luni^{-1}_g
\end{align*}
so that $E$ is a Hopf algebroid in the sense of \cite{BSz: Hgd, Day-Street}. 
Furthermore, the element $\iota\ci\pi\in E$ is a Frobenius integral in the 
ordinary sense.This explains how Theorem \ref{thm: main} implies the following

\begin{cor}
For a monoidal category $\C$ the following conditions are equivalent:
\begin{enumerate}
\item There is a Frobenius Hopf algebroid $A$ and a monoidal equivalence 
$\C\simeq\M_A$.
\item $\C$ is cocomplete, Abelian with a small projective generator  
admitting a strong Frobenius structure.
\item $\C$ is closed, cocomplete, Abelian with a small projective 
selfdual generator $g$ admitting a comonoid structure and a Frobenius 
integral and such that the map in Theorem \ref{thm: main} (4) is an isomorphism.
\end{enumerate}
\end{cor}

\section{Examples}

\subsection{Azumaya algebras}

Consider the trivial bialgebroid $E(R)$ over $R$ and the trivial bialgebroid 
over the trivial algebra $E(k)\cong k$. There is a monoidal Morita equivalence 
$E(R)\simeq E(k)$ iff $_R\M_R\simeq \M_k$ iff $R$ is an Azumaya $k$-algebra by a 
Theorem of Takeuchi \cite{Takeuchi: sqrM}. In this case an equivalence 
$P:E(R)\iso E(k)$ can be given as follows. $P$ is the $k$-module $R$ with right 
$R^e$-action $r\ract(r_1\o r_2)=r_1rr_2$, with comultiplication $\cop_P(r)=r\oR 
1\in R\oR R$ and counit $\eps_P(r)=r$. Thus $P$ is nothing else but the monoidal 
unit of $\M_{E(R)}$ with its canonical comonoid structure. Therefore it is a 
strong comonoid. Since an Azumaya algebra is a finite projective and  
split extension of $k$ \cite{Kadison}, the $k$-module $_kR$ is a progenerator. 
By centrality of $R$ the bimodule $_kP_{E(R)}$ is faithfully balanced. Therefore 
$P_{E(R)}$ is also a progenerator. This proves that the monoidal unit 
$P=\,_RR_R$ is a strong comonoid progenerator in $_R\M_R$ for any Azumaya 
algebra $R$.

\subsection{Blowing up}

Let $B$ be a bialgebra over $k$ and $n\in\NN$. It is well known that 
$A=\Mat_n(B)$ is a weak bialgebra with source and target subalgebras coinciding 
with $R:=\Diag_n(k)$. Now we consider its right bialgebroid version. Right 
multiplication of $R$ on $A$ makes it an $R$-$R$-bimodule since $R$ is 
commutative. Choosing a set of matrix units the coproduct and counit are defined 
by
\begin{align*}
\cop_A(e_{ij}\o b)&:=(e_{ij}\o b\oneT)\oR(e_{ij}\o b\twoT)\\
\eps_A(e_{ij} b)&:=\delta_{ij}\eps_B(b)
\end{align*}
where $\delta_{ij}$ is the Kronecker symbol. 

We would like to show that there is a monoidal Morita equivalence $A\simeq B$. 
For this purpose let $P=B^n$, the column vectors with entries from $B$. It is an 
$A$-$B$ bimodule via $(a\lact p\ract b)_i=\sum_j  a_{ij}p_jb$ and is a Morita 
equivalence bimodule. Next we make $P$ into a $k$-coalgebra:
\begin{align*}
\cop_P(e_i\o b)&=(e_i\o b\oneT)\o(e_i\o b\twoT)\\
\eps_P(e_i\o b)&= \eps_B(b)
\end{align*}
where $\{e_i\}$ is the canonical basis of column vectors. It is easy to verify 
that (1)-(5) of Lemma \ref{lem: 1-cell} hold, thus $P:B\to A$ is a 1-cell in 
$\BGD$. It remains to show that the comonoid $P$ is strong. This can be seen 
directly on the functor $\Ha(P)$. Indeed, $\Hom_B(P,V)\cong \oplus^n V$ with the 
obvious right $A$-module structure. But 
\begin{equation*}
\oplus^n V\ \oR\ \oplus^n W\ \cong \oplus^n(V\o W)
\end{equation*}
so $\Ha(P)$ is strong monoidal. This proves that $P_B$ is a strong comonoid 
progenerator and therefore $A\simeq B$.

\subsection{Drinfeld twists}

For a bialgebroid $B$ over $S$ a Drinfeld twist can be defined as an invertible 
element $J\in B\ex{S} B$ such that
\begin{align}
\label{twist0}
J\,\eta(s\o s')&=\eta(s\o s')\,J\qquad s,s'\in S\\
\label{twist1}
(J\oS 1)(\cop\oS B)(J)&=(1\oS J)(B\oS\cop)(J)\\
\label{twist2}
(\eps\oS B)(J)&=1\ =(B\oS\eps)(J)
\end{align}
The (\ref{twist0}) ensures that the twisted comultiplication 
$\tilde\cop(b):=J\cop(b)J^{-1}$ remains an $S$-$S$-bimodule map. It is easy to 
show that equations (\ref{twist1}) and (\ref{twist2}) imply that $\tilde\cop$ is 
coassociative and has counit $\tilde\eps=\eps$. In this way one obtains the 
twisted bialgebroid $\tilde B$ with the same underlying $S^e$-ring structure as 
$B$. 

Now we construct an equivalence $P:B\to \tilde B$. As a $\tilde B$-$B$-bimodule 
let $P=\,_BB_B$. The coring structure over $S$ is given by $\cop_P(b):=J\cop(b)$ 
and $\eps_P(b)=\eps(b)$ for $b\in P$. Coassociativity of $\cop_P$ is clear from 
(\ref{twist1}) while counitality can be shown using (\ref{twist2}) and 
(\ref{cgd2}). E.g.,
\begin{align*}
(\eps_P\oS P)\circ\cop_P(b)&=\eps(J_1b\oneT)\ud (J_2b\twoT)
=J_2b\twoT t\eps(s\eps(J_1)b\oneT)\\
&=J_2t\eps(J_1)b\twoT t\eps(b\oneT)=b\twoT t\eps(b\oneT)=b\,.
\end{align*}
The $\cop_P$ also satisfies (\ref{cgd2}) due to (\ref{twist0}) and we have
\begin{align}
\cop_P(bb')&=\tilde\cop(b)\cop_P(b')\\
\eps_P(bb')&=\eps(s\eps(b)b')=\eps_P(\tilde\eps(b)\ld b')
\end{align}
Therefore Lemma \ref{lem: 1-cell} implies that $P$ is a 1-cell in $\BGD$ from 
$B$ to $\tilde B$. Clearly, $P_B$ is a progenerator. In order to see that $P$ is 
a strong comonoid in $\M_B$, let us note at first that $T=\Hom_B(P, S)\cong  S$
acts on $E=\End_B(P)\cong B$ "like" $S$ acts on $B$. At second, we 
can use the criterion of Corollary \ref{cor: strongness} (2) as follows. 
Composing the map there with the tensor square of the isomorphism $\End 
(B_B)\cong B$ we obtain
\begin{align*}
B\oS B&\iso E\oT E\iso \Hom_B(P,P\oS P)\\
b\oS b'&\mapsto \{b''\mapsto b{b''}\oneT\oS b'{b''}\twoT\}
\end{align*}
which is an isomorphism because $P_B$ is cyclic. This proves that $P_B$ is a 
strong comonoid progenerator therefore $E\simeq B$. But also $\lambda_P:\tilde 
B\to E$ is an isomorphism of bialgebroids, hence $\tilde B\simeq B$.

Note that any pair of isocategorical groups \cite{EG: isocat groups} 
provides a special case of this example.

\subsection{\sqrM base change}

In case of trivial bialgebroids $E(R)$ and $E(S)$ our definition of monoidal 
Morita equivalence  $E(R)\simeq E(S)$ reduces to Takeuchi's \sqrM equivalence 
of $R$ and $S$ \cite{Takeuchi: sqrM}. If $B$ is a bialgebroid over $S$ and 
$P:E(R)\iso E(S)$ is a monoidal Morita equivalence bimodule then Schauenburg's 
construction in \cite[Theorem 4.3]{Schauenburg: Morita} corresponds to an 
invertible 2-cell
\begin{equation}
\begin{CD}
B@>X>>A\\
@V{\eta^*_B}V{\qquad \cong}V @VV{\eta^*_{A}}V\\
E(S)@>Q>> E(R)
\end{CD}
\end{equation}
in $\BGD$ where $Q$ is an inverse equivalence of $P$ and the bialgebroid $A$ and 
the 1-cell $X$ is constructed as follows. In the language of $\Cgd$ the 
$A$ is the monad
\begin{align}
A&:=\bra Q\odot B\odot P, \mu_{A}, \eta_{A}\ket\\
&\begin{CD}
\mu_{A}=A\odot A@>\sim>>Q\odot B\odot B\odot 
P@>Q\odot\mu_B\odot P>>A
\end{CD}\notag\\
&\begin{CD}
\eta_{A}=E(R)@>\sim>>Q\odot P@>Q\odot\eta_B\odot P>>A
\end{CD}\notag
\end{align}
and the $X$ is the bimodule
\begin{align}
X&:=\bra Q\odot B, \lambda_X,\rho_X\ket\\
&\begin{CD}
\lambda_X=A\odot X@>\sim>>Q\odot B\odot B@>Q\odot\mu_B>>X
\end{CD}\notag\\
&\begin{CD}
\rho_X=X\odot B@>Q\odot\mu_B>>X
\end{CD}\notag
\end{align}
The required invertible 2-cell $\eta_A^*\oA X\to Q\o_{S^e}\eta_B$ is obtained as 
follows. As a left $A$-module map $\lambda_X$ provides the isomorphism 
$\eta_A^*\oA X\to X$ of 1-cells in $\BGD$, on the one hand. On the other hand, 
$X$ is the horizontal composition $X=Q\o_{E(S)} \eta_B^*$ of 1-cells in $\BGD$.

The 1-cell $X$ can be easily shown to be an equivalence by constructing an 
inverse equivalence $Y:A\to B$. It is left to the reader that $Y=B\odot P$ with 
its obvious $B$-$A$-bimodule structure is indeed an inverse equivalence.

We note that bialgebroids $A$ and $B$ can be monoidally Morita equivalent 
without being there any equivalence $E(R)\simeq E(S)$ between their base
bialgebroids, i.e., without their base rings being \sqrM equivalent. A simple 
example is the blowing up of a bialgebra.


\begin{thebibliography}{XX} 
\begin{small} 

\bibitem{BSz: Hgd} G. B\"ohm, K. Szlach\'anyi, \textit{Hopf algebroids with 
bijective antipodes: axioms, integrals, and duals}, 
J. Algebra \textbf{274} (2004) 708-750

\bibitem{Brz-Mil} T.~Brzezi\'nski and G.~Militaru, \textit{Bialgebroids,
$\times_A$-bialgebras and duality}, J. Algebra \textbf{251} (2002) 279-294

\bibitem{Day-Street} B. Day, R. Street, \textit{Quantum categories, star 
autonomy, and quantum groupoids},
in "Galois theory, Hopf algebras, and semiabelian categories", Fields Inst. 
Comm. \textbf{43} (2004) 187-225

\bibitem{EG: isocat groups} P. Etingof, S. Gelaki, \textit{Isocategorical 
groups}, Internat. Math. Res. Notices 2001, no.2, 59-76

\bibitem{Kadison} L. Kadison, \textit{New examples of Frobenius extensions}, 
University Lecture Series, Vol. 14, Amer. Math. Soc., Providence, 1999

\bibitem{Leinster: BB} T. Leinster, \textit{Basic bicategories},
\texttt{math.CT/9810017}

\bibitem{Lu} J.-H. Lu, \textit{Hopf algebroids and quantum groupoids},
Int. J. Math. \textbf{7} (1996) 47-70

\bibitem{Mitchell} B. Mitchell, \textit{Theory of Categories}, Academic Press 
New York - London, 1965

\bibitem{Mueger: I} M. M\"uger, \textit{From subfactors to categories and 
topology. I. Frobenius algebras in and Morita equivalence of tensor categories}, 
J. Pure Appl. Algebra \textbf{180} (2003) 81-157

\bibitem{Pareigis: Morita} B. Pareigis, \textit{Morita equivalence of module 
categories with tensor product}, Comm. Algebra \textbf{9} (1981) 1455-1477

\bibitem{Pareigis: nature} B. Pareigis, \textit{A non-commutative
non-cocommutative Hopf algebra in "nature"}, J. Algebra \textbf{70} (1981)
356-374   

\bibitem{Schauenburg: Bial} P. Schauenburg, \textit{Bialgebras over 
noncommutative rings, and a structure theorem for Hopf bimodules}, 
Applied Categorical Structures \textbf{6} (1998) 193-222

\bibitem{Schauenburg: Duals} P. Schauenburg, \textit{Duals and Doubles of 
Quantum Groupoids}, in "New trends in Hopf algebra theory'' , Contemporary 
Mathematics \textbf{267} (2000) 273

\bibitem{Schauenburg: Morita} P. Schauenburg, \textit{Morita base change in 
quntum groupoids}, in "Locally Compact Quantum Groups and Groupoids"
ed.: L. Vainerman  (IRMA Lectures in Mathematics and Theoretical
Physics 2) de Gruyter 2003, 79-104

\bibitem{Sz: Siena} K. Szlach\'anyi, \textit{Finite quantum groupoids and 
inclusions of finite type}, Fields Inst. Comm. \textbf{30} (2001) 393-407

\bibitem{Sz: EM} K. Szlach\'anyi, \textit{The monoidal Eilenberg-Moore
construction and bialgebroids}, J. Pure Appl. Algebra \textbf{182} (2003) 
287-315 

\bibitem{Takeuchi: x} M. Takeuchi, \textit{Groups of algebras over $A \o 
\overline{A}$}, J. Math. Soc. Japan \textbf{29} (1977) 459-492  

\bibitem{Takeuchi: sqrM} M. Takeuchi, \textit{$\sqrt{\text{Morita}}$ theory - 
Formal ring laws and monoidal equivalences of categories of bimodules},
J. Math. Soc. Japan \textbf{39} (1987) 301-336


\end{small} 
\end{thebibliography}
\end{document}